\documentclass[11pt]{article}
\usepackage[a4paper,margin=28mm]{geometry}
\usepackage{needspace}
\usepackage[T1]{fontenc}
\usepackage{lmodern}
\usepackage{amsmath,amssymb,amsthm,booktabs,array,hyperref,xcolor,tikz}
\IfFileExists{microtype.sty}{\usepackage{microtype}}{}
\hypersetup{colorlinks=true,linkcolor=black,citecolor=black,urlcolor=black,pdftitle={Zero-sum Ramsey numbers of forests, cycles and complete bipartite graphs},pdfauthor={Linlin Wang and Guoqing Wang}}
\newtheorem{theorem}{Theorem}[section]
\newtheorem{lemma}[theorem]{Lemma}
\newtheorem{conjecture}[theorem]{Conjecture}
\newtheorem{remark}[theorem]{Remark}
\theoremstyle{plain}
\newcommand{\Z}{\mathbb{Z}}
\newcommand{\F}{\mathbb F}
\newcommand{\defect}{\operatorname{def}}
\newcommand{\red}{\mathrm{red}}
\newcommand{\blue}{\mathrm{blue}}
\title{Zero-sum Ramsey numbers of forests, cycles\\
and complete bipartite graphs}
\author{Linlin Wang\thanks{School of Mathematics,
China University of Mining and Technology.
Email: \href{mailto:5400@cumt.edu.cn}{\texttt{5400@cumt.edu.cn}}.}
\and Guoqing Wang\thanks{Corresponding author. School of Mathematics, Tiangong University.
Email: \href{mailto:gqwang1979@aliyun.com}{\texttt{gqwang1979@aliyun.com}}.}}
\date{}
\begin{document}
\maketitle

\begin{abstract}
For a graph $G$ and $m\ge2$ dividing $e(G)$, the zero-sum Ramsey
number $R(G,\Z_m)$ is the least $N$ such that every $\Z_m$-labeling
of $E(K_N)$ contains a copy of $G$ whose edge labels sum to zero.
We disprove the Bialostocki--Dierker asymptotic conjecture for forests
by constructing trees for every sufficiently large edge count whose
ratio of zero-sum to two-color Ramsey numbers has lower limit at
least $554/553>1$. For powers of two $q\ge2$, we characterize all
$s,k\ge1$ satisfying $R(K_{s,kq},\Z_q)=s+kq$ and prove
$R(mK_{q,q},\Z_m)=2qm+m/\gcd(m,q)-1$ for every $m\ge2$.
The packing upper bound combines Caro's transfer inequality with
our spanning result and sharpens earlier bounds when $\gcd(m,q)\ge4$.
The formula also gives exact values for counterexamples to Caro's
defect conjecture, whose failure was already implicit in earlier
work of Caro and Caro--Provstgaard. For every integer $m\ge2$ and every
multiple $\ell$ of $m$ with $\ell\ge100m^3$, we also prove
$R(C_\ell,\Z_m)=\ell+m/\gcd(m,2)-1$
for the cycle $C_\ell$ on $\ell$ vertices, determining the eventual value
left by the bounds of Chi and He for even moduli $m\ge4$.
\end{abstract}

\section{Introduction}
\label{sec:introduction}

The celebrated Erd\H{o}s--Ginzburg--Ziv theorem, proved by
P.~Erd\H{o}s, A.~Ginzburg and A.~Ziv in 1961~\cite{EGZ}, states
that every sequence of $2m-1$ elements of $\Z_m$ contains $m$ terms
whose sum is zero. Together with the Davenport constant, this theorem
motivated the development of zero-sum theory in combinatorial number
theory. The subject studies substructures with prescribed additive
properties in sequences, sets and other combinatorial structures,
mainly over finite abelian groups; see the
survey~\cite{GaoGeroldinger}, the monograph of Geroldinger and
Halter-Koch~\cite[Chapter~5]{GeroldingerHalterKoch}, and the algebraic
and counting developments in~\cite{GGSubsequenceCount,GGHGroupAlgebras}.

Zero-sum Ramsey theory extends these questions to graphs and
hypergraphs~\cite{CaroDelta,CaroProvstgaard,ZhangWang}.
For a finite simple graph $G$, write $e(G)$ for its number of edges,
and let $m\ge2$ divide $e(G)$. The zero-sum Ramsey number
$R(G,\Z_m)$ is the least $N$ such that every labeling
$w:E(K_N)\to\Z_m$ contains a copy $G'$ of $G$ with
$w(G'):=\sum_{e\in E(G')}w(e)=0$. Copies need not be induced.
For any set $F$ of labeled edges, write
$w(F)=\sum_{e\in F}w(e)$.
Write $R(G,r)$ for the ordinary $r$-color Ramsey number, the least
$N$ such that every $r$-coloring of the edges of $K_N$ contains a
monochromatic copy of $G$, and $tG$ for $t$ vertex-disjoint copies
of $G$. Taking the $m$ residues of $\Z_m$ as colors shows that
$R(G,\Z_m)\le R(G,m)$ whenever $m\mid e(G)$: a monochromatic copy
of edge label $x$ has total weight $e(G)x\equiv0\pmod m$. If
$e(G)=m$, a $\{0,1\}$-labeling shows conversely that
$R(G,2)\le R(G,\Z_m)$, since an $m$-edge copy then has weight zero
exactly when it is monochromatic. Hence (see~\cite[Theorem~3.1]{CaroSurvey})
\begin{equation}\label{eq:sandwich}
R(G,2)\le R(G,\Z_m)\le R(G,m)
\qquad\text{if }e(G)=m\ge2.
\end{equation}

Bialostocki and Dierker initiated this theory and obtained early
results for small graphs, stars and matchings~\cite{BDSmall,BD92}.
Caro obtained further results for stars and complete
graphs~\cite{CaroCliques,CaroStars}, characterized the graph case
modulo two~\cite{CaroParity}, and developed the theory of disjoint
copies and hypergraphs~\cite{CaroDelta}. We consider the following
two conjectures, recorded in Caro's survey~\cite{CaroSurvey}.

\begin{conjecture}[Bialostocki--Dierker; {\cite[Conjecture~3.1]{CaroSurvey}}]
\label{conj:forest}
\begin{equation}\label{eq:forest-conjecture}
\lim_{m\to\infty}
\max_{\substack{F\text{ a forest}\\e(F)=m}}
\frac{R(F,\Z_m)}{R(F,2)}=1.
\end{equation}
\end{conjecture}

The second conjecture uses Caro's defect $\defect_t(G)\in\{0,1\}$,
defined in Section~\ref{sec:packings}. We state its graph case.

\begin{conjecture}[Caro; {\cite[Conjecture~1]{CaroDelta}};
{\cite[Conjecture~3.4]{CaroSurvey}}]\label{conj:defect}
For every graph $G$ on $n$ vertices and every integer $t\ge2$,
\begin{equation}\label{eq:defect-conjecture}
R(tG,\Z_t)=t(n+d)-d,
\qquad d=\defect_t(G).
\end{equation}
\end{conjecture}

We disprove Conjecture~\ref{conj:forest}, characterize when certain
complete bipartite graphs attain the vertex-count bound, and obtain
exact zero-sum packing formulas and exact values for sufficiently long
cycles. These packing results sharpen
earlier bounds of Caro~\cite{CaroDelta} and Caro and
Provstgaard~\cite{CaroProvstgaard}, which already imply the failure
of Conjecture~\ref{conj:defect}; see Remark~\ref{rem:earlier-work}.

Let $T_{a,b}$ be obtained from $K_{1,a}$ by subdividing $b$ distinct
edges once, where $a\ge b\ge1$. Chi, He and
Sun~\cite[Theorem~1.6]{ChiHeSun} proved strict separation for this
family. Their theorem already gives a tree $T$ with $e(T)=m$ and
$R(T,\Z_m)>R(T,2)$ for every $m\ge57$: take $(a,b)=(m-5,5)$ for
odd $m$ and $(a,b)=(m-6,6)$ for even $m$. We strengthen this
separation to a fixed positive asymptotic ratio gap, uniformly over
all sufficiently large edge counts.

\begin{theorem}\label{thm:main-forest}
For every positive integer $k$, put $T_k=T_{276k,12k}$. Then
\[
e(T_k)=288k,\qquad
R(T_k,2)\le553k,\qquad
R(T_k,\Z_{288k})\ge554k-1.
\]
Consequently
\[
\liminf_{k\to\infty}
\frac{R(T_k,\Z_{288k})}{R(T_k,2)}\ge\frac{554}{553}>1.
\]
Moreover, for every integer $m\ge576$ there is a tree $T$ with $e(T)=m$
and $R(T,\Z_{m})>R(T,2)$, and
\[
\liminf_{m\to\infty}
\max_{\substack{F\text{ a forest}\\e(F)=m}}
\frac{R(F,\Z_{m})}{R(F,2)}\ge\frac{554}{553}>1.
\]
In particular, Conjecture~\ref{conj:forest} is false.
\end{theorem}

Chi and He~\cite[Theorem~1.1]{ChiHeBipartite} proved that
$R(K_{s,kq},\Z_q)=s+kq$ for all $k\ge1$ when $s$ is sufficiently
large. For powers of two, we determine all such side sizes and
combine this spanning result with a general inequality of Caro to
obtain an exact packing formula.

For cycles, Chi and He~\cite[Theorem~1.3 and Corollary~1.4]{ChiHeCycles}
proved that, for each fixed odd $m\ge3$ and all sufficiently large $k$,
$R(C_{mk},\Z_m)=mk+m-1$. For even $m\ge4$, their bounds leave
\[
mk+\frac m2-1\le R(C_{mk},\Z_m)\le mk+m-1
\]
when $k$ is sufficiently large. We determine the lower endpoint as
the eventual value, with an explicit length condition.

\begin{theorem}\label{thm:main-cuts}
\begin{enumerate}
\item[(i)] Let $q\ge2$ be a power of two and $s,k\ge1$. Then
\[
R(K_{s,kq},\Z_{q})=s+kq
\quad\Longleftrightarrow\quad
s\ge q\ \text{ or }\ q-s\in\{2^j:j\ge0\}.
\]
For $q\ge4$, the least $s_0$ for which equality holds for every
$k\ge1$ and every $s\ge s_0$ is $q-2$.
\item[(ii)] For every power of two $q\ge2$ and every integer $m\ge2$,
\[
R(mK_{q,q},\Z_{m})=2qm+\frac{m}{\gcd(m,q)}-1.
\]
\item[(iii)] Let $m\ge2$ be an integer, and let $\ell$ be a positive
multiple of $m$ with $\ell\ge100m^3$. Then
\[
R(C_\ell,\Z_m)=\ell+\frac{m}{\gcd(m,2)}-1.
\]
\end{enumerate}
\end{theorem}

In particular, part~(iii) gives
\[
R(C_{mk},\Z_m)=mk+\frac m2-1
\qquad(m\ge4\text{ even},\ k\ge100m^2).
\]
The length condition concerns large multiples of a fixed modulus.
Caro's conjecture $R(C_m,\Z_m)=4m-3$ for odd $m\ge5$, recorded in
\cite[footnote to Theorem~1.1]{ChiHeCycles}, concerns the separate case $\ell=m$.

\begin{remark}[Relation to earlier work]\label{rem:earlier-work}
Let $q\ge2$ be a power of two.
Caro's bounds for disjoint copies~\cite{CaroDelta}, restated in
\cite[Theorem~A]{CaroProvstgaard}, and the degree criterion in
\cite[Theorem~4]{CaroProvstgaard} already imply counterexamples to
Conjecture~\ref{conj:defect}. Indeed,
$\defect_2(K_{q,q})=0$ and $\defect_{2q}(K_{q,q})=1$;
applying Theorem~A(2) with prime factor $2$ gives
\[
R(2qK_{q,q},\Z_{2q})\le4q^2+q-1<4q^2+2q-1.
\]
The rightmost expression is the conjectured value. Moreover, when
$\gcd(m,q)\in\{1,2\}$, the exact value in
Theorem~\ref{thm:main-cuts}(ii) follows from the same earlier bounds
and the degree-count construction in the proof of
\cite[Theorem~5]{CaroProvstgaard}.

The spanning assertion for $q=2$ also follows from
\cite[Theorem~2]{CaroProvstgaard}: label each $(s+1)$-set by the
sum of the labels on its clique edges and apply that theorem to a
delta-system with an $s$-vertex core and $2k$ petals. The core edges
are counted $2k$ times, leaving the weight of a spanning $K_{s,2k}$
modulo two.

Thus the new contributions of parts~(i) and~(ii) are the full spanning
criterion for powers of two $q\ge4$ and the sharper packing values
when $\gcd(m,q)\ge4$. For the packing upper bound, we combine
Caro's transfer inequality~\cite[Theorem~8]{CaroDelta} with the
spanning value supplied by Lemma~\ref{thm:balanced-parity}; the lower
bound adapts the degree-count construction cited above.
\end{remark}

Section~\ref{sec:forests} proves Theorem~\ref{thm:main-forest}.
Section~\ref{sec:balanced-cuts} proves Theorem~\ref{thm:main-cuts},
compares the exact packing values with the defect formula, and treats
sufficiently long cycles in its final subsection.

\Needspace{24\baselineskip}
\section{Proof of Theorem~\ref{thm:main-forest}}
\label{sec:forests}

We prove Theorem~\ref{thm:main-forest} and extend its lower-bound
construction to trees with other bounded branch shapes.

\subsection{Preliminaries and the zero-sum lower bound}
\label{sec:prelim}

In a red--blue coloring, $N_{\red}(v)$ and $N_{\blue}(v)$ denote the
two neighborhoods, and $e_{\red}(A,B)$ counts red edges between
disjoint sets $A$ and $B$.

For forests with $3\mid e(F)$, Caro and
Mifsud~\cite[Corollary~11]{CaroMifsud} proved
$R(F,\Z_3)\le |V(F)|+2$, and Alvarado, Colucci and
Parente~\cite[Theorem~1.1]{AlvaradoColucciParente} determined the exact
values without isolated vertices.

Related tree questions were studied by Caro and
Roditty~\cite{CaroRoditty}. Linear bounds for forests at prime
moduli~\cite{ColucciDEmidio} and for bounded-degree
graphs~\cite{KatzEtAl} give further context.
For subdivided stars, Zhao, Mao and
Zhou~\cite[Theorem~1.6]{ZhaoMaoZhou} determine exact ordinary
Ramsey numbers under conditions requiring $a$ to grow at least
quadratically with $b$. Our estimates apply when $a$ and $b$
grow proportionally.

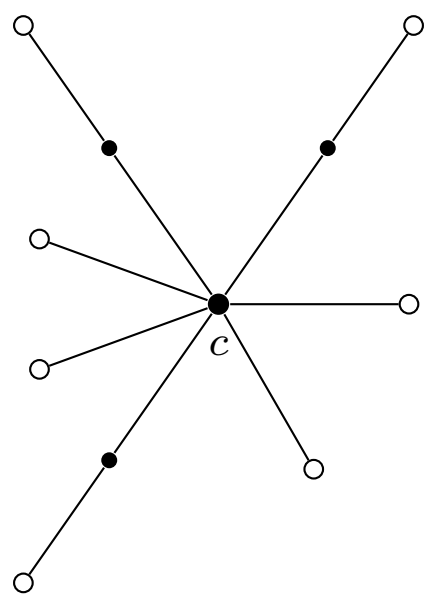
\begin{figure}[htbp]
\centering
\begin{tikzpicture}[scale=0.9,
  ctr/.style={circle,fill=black,inner sep=0pt,minimum size=4pt},
  mid/.style={circle,fill=black,inner sep=0pt,minimum size=3pt},
  lf/.style={circle,draw,fill=white,inner sep=0pt,minimum size=3.5pt}]
  \node[ctr,label=below:{$\scriptstyle c$}] (c) at (0,0) {};
  \foreach \a in {55,125,235}{
    \node[mid] (m\a) at (\a:1.4) {};
    \node[lf]  (l\a) at (\a:2.5) {};
    \draw (c)--(m\a)--(l\a);
  }
  \foreach \a in {0,160,200,300}{
    \node[lf] (d\a) at (\a:1.4) {};
    \draw (c)--(d\a);
  }
\end{tikzpicture}
\caption{The tree $T_{a,b}$. The center $c$ is adjacent to $a$
vertices; $b$ of them are subdivision vertices, drawn as filled dots,
each carrying one further leaf. Three subdivided arms and four direct
arms are shown. The $b$ terminal edges are pairwise disjoint.}
\label{fig:tab}
\end{figure}

The $b$ subdivided arms of $T_{a,b}$ have length two; the remaining
$a-b$ arms have length one. Its bipartition classes have orders
$a$ and $b+1$, and apart from the center $c$ every vertex has degree at
most two; see Figure~\ref{fig:tab}.

\paragraph{The lower-bound construction.}

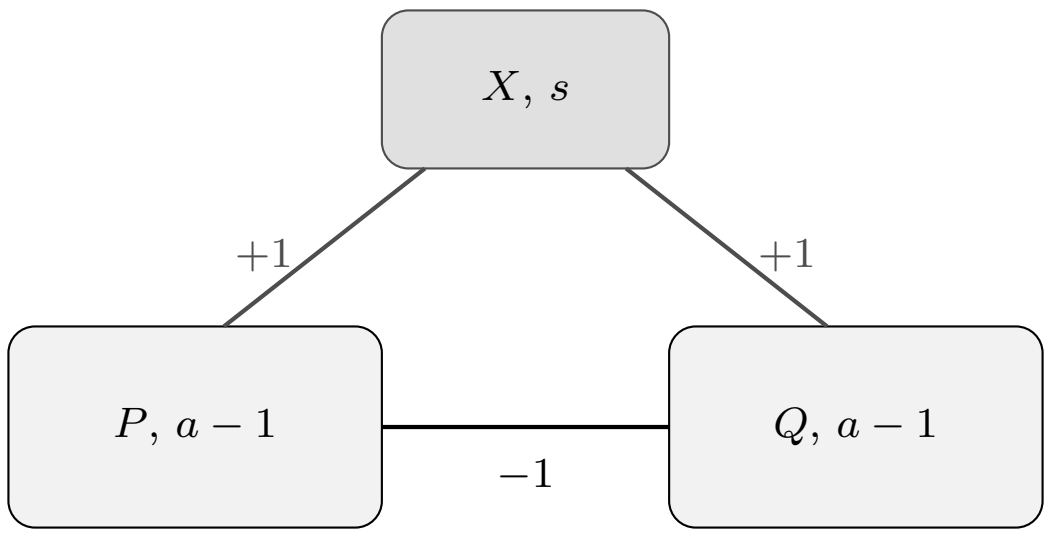
\begin{figure}[htbp]
\centering
\begin{tikzpicture}[scale=0.95,every node/.style={font=\scriptsize}]
  \draw[rounded corners=5pt,fill=gray!10] (0,0) rectangle (2.6,1.4);
  \node at (1.3,0.7) {$P$, $a-1$};
  \draw[rounded corners=5pt,fill=gray!10] (4.6,0) rectangle (7.2,1.4);
  \node at (5.9,0.7) {$Q$, $a-1$};
  \draw[rounded corners=5pt,fill=black!12,draw=black!70]
    (2.6,2.5) rectangle (4.6,3.6);
  \node at (3.6,3.05) {$X$, $s$};
  \draw[thick] (2.6,0.7)--(4.6,0.7) node[midway,below=2pt]{$-1$};
  \draw[thick,black!70] (1.5,1.4)--(2.9,2.5) node[pos=0.45,left]{$+1$};
  \draw[thick,black!70] (5.7,1.4)--(4.3,2.5) node[pos=0.45,right]{$+1$};
\end{tikzpicture}
\caption{The labeling of Lemma~\ref{lem:forest-lower}. Edges meeting
$X$ carry $1$, edges between $P$ and $Q$ carry $-1$, and all remaining
edges carry $0$.}
\label{fig:labeling}
\end{figure}

The construction below extends the labeling of Chi, He and
Sun~\cite[Proposition~5.1]{ChiHeSun} by replacing the single
exceptional vertex with a set $X$; see Figure~\ref{fig:labeling}.

\Needspace{19\baselineskip}
\begin{lemma}\label{lem:forest-lower}
Let $a>b\ge1$ and $s\ge1$ be integers.
\begin{enumerate}
\item[(i)] If $6s\le b+1$, then
\[
R(T_{a,b},\Z_{a+b})\ge2a+s-1.
\]
More generally, let $T$ be a tree with $m$ edges and bipartition sizes
$t_1\le t_2$, let $L,s\ge1$, and suppose that for some vertex $c$ each
component of $T-c$ has at most $L$ vertices and that $T-c$ contains a
matching of size $s$. If
\[
\deg_T(c)>m/2 \qquad\text{and}\qquad L(L+1)s\le t_1,
\]
then $R(T,\Z_m)\ge2t_2+s-1$.
\item[(ii)] Conversely, let $a\ge4b$ and $a,b,s\ge1$. For the labeling
used in (i), with $|P|=|Q|=a-1$ and $|X|=s$, there is no zero-sum copy
of $T_{a,b}$ if and only if $6s\le b+1$.
\end{enumerate}
\end{lemma}

\begin{proof}
(i) We prove the general bound first. Use the labeling in
Figure~\ref{fig:labeling} with $|P|=|Q|=t_2-1$ and $|X|=s$.
Edges meeting $X$ have label $1$, edges between $P$ and $Q$ have
label $-1$, and the other edges have label $0$.

For a proposed copy let $h$ count its edges meeting $X$, and let $r$
count its edges between $P$ and $Q$. Its weight is $h-r$ modulo $m$.
If the copy avoids $X$, zero weight requires $r=0$ or $r=m$. In the
first case one part contains the whole tree. In the second one part
contains a bipartition class of size $t_2$. Both exceed the part size $t_2-1$.

If $c$ lies in $X$, then $h\ge\deg_T(c)>m/2$. At most $s-1$ other
vertices lie in $X$, so one of the $s$ disjoint edges in $T-c$ misses
$X$. Hence $h<m$. Since $r\le m-h<h$, we have $0<h-r<m$.

It remains to consider a copy that meets $X$ but has $c$ outside $X$.
Every vertex other than $c$ has degree at most $L$, so $1\le h\le Ls<m$.
As $r\le m-h<m$, zero weight forces $r=h$. Discard each component
of $T-c$ that meets $X$ or that contains an endpoint other than $c$
of an edge between $P$ and $Q$. At most $s+r$ components are discarded,
since the endpoints of each such edge other than $c$ lie in one component.
The remaining subtree lies in one part, yet has at least
\[
t_1+t_2-L(s+r)\ge t_1+t_2-L(L+1)s\ge t_2
\]
vertices. Thus the labeling avoids a zero-sum $T$ on $2t_2+s-2$ vertices.

For $T_{a,b}$, take $L=2$, $t_1=b+1$, and $t_2=a$.
The center has degree $a>(a+b)/2$, and the $b$ terminal edges give
the required matching since $s\le b$. The condition $L(L+1)s\le t_1$
is exactly $6s\le b+1$. The labeling uses $2a+s-2$ vertices, as required.

(ii) Part (i) proves sufficiency. Now suppose
$6s\ge b+2$, and put $u=\lceil(b+2)/6\rceil\le\min\{s,b\}$.
Place the center in $P$. Send $u$ subdivided arms through $X$ to $Q$.
Their edges have total weight $2u$.

To cancel this weight, place $2u$ further arms entirely in $Q$
apart from their common center. Each contributes $-1$, regardless of
its length. There are enough arms since $a-u\ge2u$. Choose subdivided
arms first. Then exactly $v=\min\{2u,b-u\}$ of the chosen arms have
a second vertex in $Q$. Put every remaining vertex in $P$.

It remains to check the part sizes. We used $u\le s$ vertices of $X$
and $3u+v\le b+2u\le3b\le a-1$ vertices of $Q$. In total we placed
\[
4u+v=\min\{6u,b+3u\}\ge b+2
\]
vertices outside $P$, so at most $(a+b+1)-(b+2)=a-1$ lie in $P$.
The placement gives the required zero-sum copy.
\end{proof}

Part (ii) of the lemma determines the exact scope of the labeling, rather
than the exact zero-sum Ramsey number. A larger exceptional set cannot
improve the lower bound without changing the construction.

\subsection{An ordinary Ramsey upper bound}

\begin{lemma}[K\"onig's theorem \cite{Konig}]\label{lem:konig}
In a finite bipartite graph, the maximum size of a matching equals
the minimum size of a vertex cover.
\end{lemma}

The following estimate refines the matching and vertex-cover
argument of Chi, He and Sun~\cite[Lemma~5.2 and Proposition~5.3]{ChiHeSun},
allowing $a$ and $b$ to grow proportionally.

\begin{lemma}\label{lem:forest-upper}
Let $a,b,q$ be positive integers such that $a\ge3b$ and
$(q+2)(a+b-2)>2b(b-1)$. Then
\[
R(T_{a,b},2)\le2a+q.
\]
\end{lemma}

\begin{proof}
Set $N=2a+q$, and suppose that a red--blue coloring of $K_N$ contains
no monochromatic $T=T_{a,b}$. Since $N-1\ge2a$, after interchanging
the colors if necessary, we may choose a vertex $v$ with at least
$a$ red neighbors. Put $A=N_{\red}(v)$.

\smallskip
\noindent\emph{A large red neighborhood.}
First suppose that $|A|\ge a+b$. Consider the red edges in $K_N-v$
that have at least one endpoint in $A$. A matching of size $b$ among
these edges gives $b$ red arms of length two centered at $v$: for each
matching edge choose an endpoint in $A$ as its middle vertex. At most
$2b$ vertices of $A$ are used, leaving at least $a-b$ direct red
leaves. This would give a red $T$.

\Needspace{9\baselineskip}
There is therefore no such matching. The endpoints $S$ of a maximal
matching cover all these red edges and satisfy $|S|\le2b-2$.
Every edge from $A\setminus S$ to another vertex outside
$\{v\}\cup S$ is blue. We have
\[
|A\setminus S|\ge a-b+2\ge b+1
\]
and
\[
N-1-|S|-(b+1)\ge2a+q-3b\ge a.
\]
Choose a center and $b$ outer endpoints in $A\setminus S$.
Then choose $a$ additional vertices outside $\{v\}\cup S$ and those
already chosen. Use $b$ of these additional vertices as the middle
vertices of the subdivided arms and the remaining $a-b$ as direct
leaves. Every required edge joins a vertex of $A\setminus S$ to a
vertex outside $\{v\}\cup S$, and hence is blue. This gives a blue
copy of $T$, again a contradiction.
This argument applies to every vertex with at least $a+b$ red
neighbors. Thus no such vertex can occur in the assumed coloring.

\smallskip
\noindent\emph{A complete blue bipartite subgraph.}
We may now write $|A|=a+j$, where $0\le j\le b-1$. Put
$B=N_{\blue}(v)$. A red matching of size $b$ between $A$ and $B$
would extend to a red $T$ centered at $v$. By Lemma~\ref{lem:konig}, the red bipartite graph between $A$ and $B$ consequently has
a vertex cover $C\subseteq A\cup B$ of order $p\le b-1$. Let
\[
X=A\setminus C,\qquad Y=B\setminus C.
\]
All edges between $X$ and $Y$ are blue, and
\begin{equation}\label{eq:XY-lower}
|X|\ge a-b+1,\qquad |Y|\ge a+q-2b+1.
\end{equation}

\smallskip
\noindent\emph{Where a blue center could lie.}
Every vertex of $X\cup Y$ has blue degree at most $a-1$. To see this,
let $U$ be the side containing a proposed center $u$ and $V$ the
other side. By \eqref{eq:XY-lower}, both sides have at least $b+1$
vertices. Also
\[
|U|+|V|=2a+q-1-p\ge a+b+1.
\]
If $u$ has at least $a$ blue neighbors, select $a$ of them, taking as
many as possible from $V$. At least $b$ selected vertices lie in $V$.
There are also $b$ unused vertices in $U$ besides $u$. If
$|V|\ge a$, all of $U\setminus\{u\}$ is available. If $|V|<a$, at
least $|U|+|V|-a-1\ge b$ are available. Pair $b$ selected vertices
in $V$ with $b$ unused vertices in $U$ to extend $b$ arms. The other
$a-b$ selected neighbors are direct leaves, giving a blue $T$.

Every $y\in Y$ already has the blue neighbors $X\cup\{v\}$, and
every $x\in X$ has the blue neighbors $Y$. The blue degree bound
therefore gives
\begin{equation}\label{eq:deficits}
|X|=a-1-\tau,\qquad |Y|=a-h,\qquad
\tau,h\ge1,\qquad p=\tau+h+q\le b-1.
\end{equation}
The identity for $p$ counts the vertices in $\{v\},C,X,Y$.

\smallskip
\noindent\emph{Counting red edges into the cover.}
Each $x\in X$ has at most $h-1$ blue neighbors in $C$, and
each $y\in Y$ has at most $\tau-1$. They therefore have at least
$p-h+1=\tau+q+1$ and $p-\tau+1=h+q+1$ red neighbors in $C$,
respectively. Summing gives
\begin{align}
e_{\red}(C,X\cup Y)
&\ge(\tau+q+1)(a-1-\tau)+(h+q+1)(a-h)\notag\\
&=(p+q+2)a-L,\label{eq:red-incidences}
\end{align}
where
\[
L=(\tau+q+1)(\tau+1)+(h+q+1)h.
\]
Using $p=\tau+h+q$ and $\tau,h\ge1$, we obtain the sharper estimate
\[
L=p(p-q)+q+2-(h-1)(2\tau-1)\le p(p-q)+q+2.
\]
If every vertex in $C$ had at most $a+b-1$ red neighbors in $X\cup Y$,
then \eqref{eq:red-incidences} would imply
\begin{align*}
(q+2)a
&\le L+p(b-1)\\
&\le p(p+b-q-1)+q+2\\
&\le2(b-1)^2-q(b-2)+2.
\end{align*}
Here the last inequality uses $p\le b-1$ and $p-q\ge2$.
Adding $(q+2)(b-2)$ gives
\[
(q+2)(a+b-2)\le2b(b-1),
\]
contrary to the hypothesis. Some vertex of $C$ therefore has at least
$a+b$ red neighbors. The large-neighborhood argument at the start of
the proof gives the required contradiction.
\end{proof}

\subsection{Asymptotic separation}

Lemmas~\ref{lem:forest-lower} and~\ref{lem:forest-upper} give
\[
R(T_{a,b},\Z_{a+b})\ge2a+s-1
\quad\text{and}\quad
R(T_{a,b},2)\le2a+q
\]
under their respective hypotheses. Subtracting the second bound from
the first gives $s-1-q$, a lower bound on the Ramsey-number gap.
For fixed $a/b=\lambda$, choosing $s$ and $q$ near the permitted
bounds gives $s-1-q=(1/6-2/(\lambda+1))b+O(1)$, which is positive
for all sufficiently large $b$ when $\lambda>11$.

\begin{lemma}\label{cor:proportional}
Let $a,b,s,q$ be positive integers such that
\[
6s\le b+1,\qquad a\ge3b,\qquad
(q+2)(a+b-2)>2b(b-1),\qquad s\ge q+2.
\]
Then $R(T_{a,b},\Z_{a+b})>R(T_{a,b},2)$.
\end{lemma}
\begin{proof}
The two lemmas give
\[
R(T_{a,b},\Z_{a+b})\ge2a+s-1\ge2a+q+1>2a+q\ge R(T_{a,b},2).
\]
\end{proof}

\begin{lemma}[Every sufficiently large edge count]\label{cor:all-moduli}
For $m\ge288$, write $m=288k+r$, where $k\ge1$ and $0\le r<288$,
and put $U_m=T_{276k+r,12k}$. Then $e(U_m)=m$ and
\[
R(U_m,2)\le553k+2r,\qquad
R(U_m,\Z_m)\ge554k+2r-1.
\]
In particular $R(U_m,\Z_m)>R(U_m,2)$ for every $m\ge576$, and
\[
\liminf_{m\to\infty}\frac{R(U_m,\Z_m)}{R(U_m,2)}
\ge\frac{554}{553}.
\]
\end{lemma}
\begin{proof}
Set $a=276k+r$, $b=12k$, $s=2k$, and $q=k$. Then
$6s\le b+1$ and $a\ge3b$. Moreover,
\[
(q+2)(a+b-2)-2b(b-1)=598k-4+(k+2)r>0.
\]
The two lemmas give the stated estimates. Their difference is $k-1$,
which is positive for $m\ge576$. Finally,
\[
\frac{R(U_m,\Z_m)}{R(U_m,2)}
\ge\frac{554k+2r-1}{553k+2r}
\longrightarrow\frac{554}{553},
\]
since $r$ remains bounded as $m\to\infty$.
\end{proof}

\paragraph{The range of proportional parameters.}
Fix a rational number $\lambda>11$ and take $a=\lambda b$, with
$b\to\infty$ through values making $a$ integral. Set
\[
s=\left\lfloor\frac{b+1}{6}\right\rfloor,
\qquad
q=\max\left\{1,
\left\lfloor\frac{2b(b-1)}{a+b-2}\right\rfloor-1\right\}.
\]
For sufficiently large $b$, these satisfy the hypotheses of both
lemmas. The floor in the definition of $q$ ensures the strict
inequality required by Lemma~\ref{lem:forest-upper}. Since
$s/b\to1/6$ and $q/b\to2/(\lambda+1)$, the lower bound for the ratio
tends to
\[
\frac{2\lambda+1/6}{2\lambda+2/(\lambda+1)}
=1+\frac{\lambda-11}{12(\lambda^2+\lambda+1)}>1.
\]
The integral choice $\lambda=23$ and $b=12k$, with $s=2k$ and $q=k$,
gives $554/553$. We choose these parameters for a simple statement;
the displayed expression is maximized over real $\lambda>11$ at
$\lambda=11+\sqrt{133}$. Part (ii) of Lemma~\ref{lem:forest-lower} shows
that any further improvement in the permitted size of $s$ requires
a different labeling. The sufficient ordinary Ramsey estimate may
still admit improvement, so the threshold $11$ is not asserted to
be an exact transition for the Ramsey numbers themselves.

\paragraph{The order of the ordinary Ramsey excess.}
Das, Reed and Skokan use subdivided stars in the proof of
\cite[Proposition~2.1]{DasReedSkokan}; specializing their two-clique
construction \cite[Lemmas~2.6(b), 2.7 and 2.8]{DasReedSkokan} with
integer parameters shows that, for $a\ge3b$ and $b\ge4$, putting
$h=\lfloor b/4\rfloor$ and $d=\lfloor h^2/(a+h)\rfloor$ gives
\[
R(T_{a,b},2)\ge2a+2d-1 .
\]
For this specialization, take
$\alpha=h$, $\beta=a+d-1-h$, and $\Delta=h-d$.
The deterministic construction in Lemma~2.6(b) of~\cite{DasReedSkokan}
applies to these integer parameters. Since
$h^2\ge d(a+d-1)$, we have
\[
 \left\lfloor\frac{\alpha\Delta}{\beta}\right\rfloor\ge d,
 \qquad
 \alpha+\beta+2(\alpha+\Delta)=a+4h-d-1\le a+b.
\]
Thus their coloring on $2a+2d-2$ vertices has no monochromatic
$T_{a,b}$.
Hence $R(T_{a,b},2)-2a=\Theta(b^2/a)$ along every sequence with
$a\ge3b$ and $b^2/a\to\infty$, and this order is attained up to the
constant by Lemma~\ref{lem:forest-upper}. For $\sqrt a=o(b)$ and
$b=o(a^{2/3})$ the excess is $o(\sqrt b)$, whereas for the same
bipartition sizes the largest excess among all trees is $\Theta(\sqrt b)$
by \cite[Theorem~1.1]{DasReedSkokan}.

\Needspace{9\baselineskip}
\begin{proof}[Proof of Theorem~\ref{thm:main-forest}]
Take $r=0$ in Lemma~\ref{cor:all-moduli} to obtain the bounds for
$T_k$. For every $m\ge288$, the maximum in
\eqref{eq:forest-conjecture} is at least the ratio for $U_m$.
The lemma therefore gives the asserted uniform lower limit,
which is greater than one and disproves Conjecture~\ref{conj:forest}.
\end{proof}

\section{Proof of Theorem~\ref{thm:main-cuts}}
\label{sec:balanced-cuts}

Caro~\cite{CaroBipartite} and Caro and
Yuster~\cite{CaroYusterBipartite} studied zero-sum bipartite Ramsey
numbers with complete bipartite hosts. Here the hosts are complete
graphs. Algebraic and counting methods have been developed in
\cite{AlonCaro,CaroLinear,CaroBinomial}, and topological methods
provide another approach~\cite{FrickEtAl}.
We first prove the spanning and packing assertions by counting cut
weights. The last subsection treats sufficiently long cycles by a
switching argument.

\subsection{Parity of cut weights}

A spanning copy of $K_{q,q}$ in $K_{2q}$ is determined by an
unordered partition into two $q$-element sets. When $q$ is a power
of two, we can determine which finite abelian groups force such a
cut to have weight zero. We count the choices of one side inside a
prescribed set of vertices. This also gives spanning complete
bipartite graphs with unequal sides and exact packing thresholds.

For a vertex set $A$ in a graph $G$, write $\delta_G(A)$ for the set
of edges with exactly one endpoint in $A$.
For a finite abelian group $\Gamma$, the Davenport constant $D(\Gamma)$
is the least positive integer such that every sequence of $D(\Gamma)$
elements of $\Gamma$ has a nonempty subsequence with sum zero.
For a finite abelian $2$-group
$\Gamma=\bigoplus_{i=1}^s\Z_{2^{a_i}}$, Olson's theorem~\cite{Olson}
gives
\[
  D(\Gamma)=1+\sum_{i=1}^s(2^{a_i}-1).
\]
We allow the empty decomposition, so $D(\Gamma)=1$ for the trivial
group.

\paragraph{Counting cuts modulo two.}
For a family $\mathcal A$ of candidate sides and $E_0\subseteq E(G)$,
let $N(E_0)$ count the sets $A\in\mathcal A$ such that
$E_0\subseteq\delta_G(A)$. We use the following criterion. If
$|\mathcal A|$ is odd and $N(E_0)$ is even whenever
$1\le |E_0|<D(\Gamma)$, then an odd number of sides in $\mathcal A$ give
weight zero and an even number give each other weight.

To see this, work in the group algebra $\F_2[\Gamma]$, with basis
elements $[g]$ satisfying $[g][h]=[g+h]$ and $[0]=1$. For a cut
$C$, expand $[w(C)]=\prod_{e\in C}(1+([w(e)]-1))$. Summing over
the candidate sides and collecting the term for each edge set gives
\begin{equation}\label{eq:cut-parity-expansion}
 \sum_{A\in\mathcal A}[w(\delta_G(A))]
   =\sum_{E_0\subseteq E(G)}N(E_0)\prod_{e\in E_0}([w(e)]-1).
\end{equation}
On the left, the coefficient of $[g]$ is the parity we seek. On
the right, the empty set contributes $[0]$, and the assumed parity
of $N(E_0)$ removes every term with $1\le |E_0|<D(\Gamma)$.

All larger products vanish in $\F_2[\Gamma]$. For a generator
$g_i$ of the $i$th cyclic summand put $y_i=[g_i]-1$. Then
$y_i^{2^{a_i}}=[2^{a_i}g_i]-1=0$ in characteristic two.
Since every $[g]$ is a product of powers of $1+y_i$, each
$[g]-1$ is a polynomial in the $y_i$ with zero constant term.
In a product of at least $D(\Gamma)=1+\sum_i(2^{a_i}-1)$ such factors,
every monomial contains some $y_i^{2^{a_i}}$. For the trivial
group, each factor is already zero. This proves the criterion.

\begin{lemma}\label{lem:prescribed-cut}
Let $q\ge2$ be a power of two and $\Gamma$ a finite abelian $2$-group
with $D(\Gamma)\le q$. Let $t$ be a positive multiple of $q$, and let
$D(\Gamma)\le h\le q$. For any $\Gamma$-labeling $w$ of the edges of a graph
$G$ and any $W\subseteq V(G)$ of size $t+h-1$, an odd number of
$t$-subsets $A$ of $W$ satisfy $w(\delta_G(A))=0$, and an even
number give each other weight. The minimum size $t+D(\Gamma)-1$ of $W$
in this assertion is sharp for guaranteeing a zero-weight cut with
the chosen $t$-set contained in $W$.
\end{lemma}
\begin{proof}
Take $\mathcal A$ to be the $t$-subsets of $W$. For $A\in\mathcal A$
write $x_i=1$ if $i\in A$ and $x_i=0$ otherwise. Over $\F_2$,
$x_i+x_j$ records whether $ij$ crosses the cut. Thus $N(E_0)$ modulo
two is the sum over $\mathcal A$ of
$P_{E_0}=\prod_{ij\in E_0}(x_i+x_j)$. Set $x_i=0$ outside $W$ and use
$x_i^2=x_i$ to make each monomial squarefree.

A monomial of degree $j$ counts the choices of $A$ containing its
$j$ specified vertices. There are $\binom{t+h-1-j}{h-1}$ choices,
since we must omit $h-1$ of the remaining vertices. When
$0\le j<D(\Gamma)\le h$, the parity of this number is the coefficient
of $X^{h-1}$ in
\[
 (1+X)^{t+h-1-j}=(1+X^q)^{t/q}(1+X)^{h-1-j}.
\]
Since $h-1<q$, only the constant term of the first factor matters.
The coefficient is therefore $1$ for $j=0$ and $0$ for $j>0$.
In particular $|\mathcal A|$ is odd. For $1\le |E_0|<D(\Gamma)$ every
monomial of $P_{E_0}$ has degree between $1$ and $|E_0|$, so $N(E_0)$ is
even. The criterion now gives the claimed parities.

For sharpness when $D(\Gamma)\ge2$, take a set $W$ of size $t+D(\Gamma)-2$ and
a vertex $v$ outside it. Label the edges from $v$ to $W$ by
$t-1$ zeros together with $2^{a_i}-1$ copies of each generator
$g_i$, and label all remaining edges by zero. Every $t$-set in
$W$ selects a nonempty submultiset of these generators. Its sum
is nonzero, since each generator occurs fewer times than its order.
When $D(\Gamma)=1$, at least $t$ eligible vertices are necessary to choose
a $t$-set.
\end{proof}

\subsection{Spanning complete bipartite graphs}

For elementary abelian groups, Alvarado, Colucci, Parente and
Souza~\cite{AlvaradoColucciParenteSouza} determine $R(C_4,\Z_2^d)$
for $d\ge2$ and study complete bipartite graphs as $d$ grows.
Here the graph has the same order as the host, and the group may
have cyclic summands of different orders.

\begin{lemma}\label{thm:balanced-parity}
Let $q\ge2$ be a power of two, and let $\Gamma$ be a finite abelian group.
Every $\Gamma$-labeling of $E(K_{2q})$ contains a zero-sum $K_{q,q}$ if
and only if $\Gamma$ is a $2$-group with $D(\Gamma)\le q$. Whenever this
condition holds, the number of unordered balanced cuts of weight
zero is odd, and the number of each other weight is even.
\end{lemma}
\begin{proof}
Fix a vertex $v$. Each unordered balanced cut has a unique side of
size $q$ avoiding $v$, so the positive assertions follow by taking
$t=h=q$ and $W=V(K_{2q})\setminus\{v\}$ in
Lemma~\ref{lem:prescribed-cut}.

Conversely, if $\Gamma$ has an element of odd prime order, label
every edge by that element. Every cut has $q^2$ edges and hence
nonzero weight. It remains to consider a $2$-group with $D(\Gamma)>q$.
Let $g_i$ generate its $i$th cyclic summand. The multiset
consisting of $2^{a_i}-1$ copies of each $g_i$ has
size at least $q$, and every nonempty submultiset has nonzero sum.
Choose $q$ of its entries and append $q-1$ zeros. Assign these
$2q-1$ labels to the edges at $v$, and give all other edges label
zero. A $q$-set avoiding $v$ selects exactly $q$ labels at $v$,
including at least one nonzero label. Its cut therefore has
nonzero weight.
\end{proof}

\begin{lemma}\label{thm:rectangular-cuts}
Let $q\ge2$ be a power of two and let $\Gamma$ be a finite abelian
$2$-group with $D(\Gamma)\le q$. If $k,s\ge1$ and $s\ge D(\Gamma)-1$, then
every $\Gamma$-labeling of $K_{kq+s}$ contains a spanning zero-sum copy
of $K_{kq,s}$. The same conclusion holds when $1\le s<q$ and
$q-s$ is a power of two, including $1$.

For $\Gamma=\Z_q$ and all positive integers $k,s$, one has
\[
 R(K_{kq,s},\Z_q)=kq+s\quad\Longleftrightarrow\quad
 s\ge q\ \text{ or }\ q-s\in\{2^j:0\le j<a\},\qquad q=2^a .
\]
In particular, for $q\ge4$ the least positive $s_0$ for which
equality holds for all $k\ge1$ and $s\ge s_0$ is $q-2$.
\end{lemma}
\begin{proof}
\emph{Existence.}
If $s\ge D(\Gamma)-1$, the host contains a set $W$ of $kq+D(\Gamma)-1$
vertices. Lemma~\ref{lem:prescribed-cut}, with $t=kq$ and
$h=D(\Gamma)$, gives a zero-weight cut whose $kq$-vertex side lies in $W$.

Now let $L=q-s$ be a power of two and consider all $s$-sets in the
host. By the binary parity rule, $\binom{M}{d}$ is odd exactly
when every $1$-bit of $d$ is also a $1$-bit of $M$. In
$\F_2[X]$, each $1$-bit $2^i$ of $M$ contributes one factor
$1+X^{2^i}$ to $(1+X)^M$. In particular, $\binom{kq+s}{s}$ is odd.

Fix a vertex $v$. For a chosen $s$-set, write $y_i=1$ if $i$ lies
on the opposite side from $v$, and $y_i=0$ otherwise. The crossing
indicator of $ij$ is still $y_i+y_j$ over $\F_2$. Thus for
$1\le |E_0|<q$ the indicator that every edge of $E_0$ crosses expands
into monomials $\prod_{i\in B}y_i$, with $1\le |B|\le |E_0|$.
Here $y_v=0$ and repeated factors satisfy $y_i^2=y_i$. It suffices
to show that each fixed nonempty $B$ lies opposite $v$ an even
number of times.

Write $j=|B|$ and $M=kq+s-j-1$. After fixing $v$ and $B$ on
opposite sides, there are $M$ free vertices. The side containing
$v$ needs either $s-1$ or $kq-1$ more vertices. Hence the count is
\[
 \binom{M}{s-1}+\binom{M}{kq-1}.
\]
Each summand is odd exactly when $j=s$. To see this, observe that
$s-1=(q-1)-L$. Its binary expansion has all the $1$-bits of $q-1$
except the bit $L$. The parity rule therefore makes the first
summand odd exactly when $M\bmod q$ is $s-1$ or $q-1$. Since
$M\equiv s-1-j\pmod q$ and $1\le j<q$, this happens only at
$j=s$. The second summand can be odd only when $M\bmod q=q-1$,
because $kq-1$ has all the low bits of $q-1$. At $j=s$ we have
$M=kq-1$, so both summands are odd.
The count is always even. Since $D(\Gamma)\le q$, the criterion
\eqref{eq:cut-parity-expansion} now gives a zero-weight cut.

\emph{The obstruction for $\Z_q$.}
Suppose $1\le s<q$ and $L=q-s$ is not a power of two. Mark $q-1$
vertices, and let $b$ be the number of them on the $kq$-vertex side
of a cut. The $s$-vertex side can hold at most $s$ marks, so
$L-1\le b\le q-1$. We will make the cut weight nonzero throughout
this interval. First label each edge by its number of marked
endpoints. Counting the edges at each marked vertex gives weight
\[
 sb+kq(q-1-b)\equiv sb\pmod q .
\]
If $s$ is odd, it is invertible modulo $q$, so this weight is
nonzero for every allowed $b$.

For even $s$, put $g=\gcd(s,q)=\gcd(L,q)$. Then $g$ is even and
$L>g$, since $L$ is not a power of two. Scale the preceding
labeling by a solution $\alpha$ of $s\alpha\equiv g\pmod q$,
then add $1$ to every edge with two marked endpoints. Exactly
$b(q-1-b)$ such edges cross, giving cut weight
\[
 gb+b(q-1-b)\equiv b(g-1-b)\pmod q .
\]
The two factors have opposite parity. Hence their product can
vanish modulo the power of two $q$ only if $b\equiv0$ or
$b\equiv g-1\pmod q$. Both residues lie below $L-1$.
This proves the equivalence. For $q\ge4$, $s=q-3$ fails the
criterion and every $s\ge q-2$ satisfies it, giving the stated
threshold.
\end{proof}

For arbitrary cyclic modulus $q$, Chi and He
\cite[Theorem~1.1]{ChiHeBipartite} prove the same spanning
conclusion for $s\ge S(q)$; for even $q\ge4$ their bound is
$S(q)=2R(K_q,\Z_{q/2})-q-1\ge q+1$, while
Lemma~\ref{thm:rectangular-cuts} gives the exact threshold $q-2$ for
powers of two $q\ge4$. Chi and He also determine the values for
$K_{s,3k}$ with $s\ge2$, including $R(K_{3,3},\Z_3)=7$
\cite[Theorem~1.2]{ChiHeBipartite}.

The next lemma shows why the hypothesis that $q$ is a power of two
cannot be removed.

\begin{lemma}\label{lem:odd-primes}
For every odd prime $p$, one has $R(K_{p,p},\Z_p)\ge2p+1$.
\end{lemma}
\begin{proof}
Assign to the $2p$ vertices the vectors
\[
e_1\ (p-1\text{ times}),\quad e_2\ (p-1\text{ times}),\quad
e_1+e_2,\quad 0
\]
in $\F_p^2$, and write these vectors as $(x_i,y_i)$. Give the edge
$ij$ weight $x_ix_j-\nu y_iy_j$, where $\nu$ is a nonsquare in
$\F_p$. For a balanced cut with sides $A$ and $B$, the total sums of
both coordinates over all vertices are zero. Its weight is therefore
\[
-\Big(\sum_{i\in A}x_i\Big)^2+
\nu\Big(\sum_{i\in A}y_i\Big)^2.
\]
Since $\nu$ is a nonsquare, this can vanish only if both coordinate
sums are zero. If $A$ omits the vertex labeled $e_1+e_2$, each
coordinate count is at most $p-1$, so both counts must be zero;
then $A$ contains at most the vertex labeled $0$. If $A$ contains
the vertex labeled $e_1+e_2$, both counts must equal $p$, forcing
$A$ to contain all $2p-2$ vertices labeled $e_1$ or $e_2$ as well.
Both alternatives contradict $|A|=p$. Thus no balanced cut has
weight zero.
\end{proof}

\subsection{Zero-sum packings and the defect formula}
\label{sec:packings}

Ordinary Ramsey numbers for disjoint copies were studied by Burr,
Erd\H{o}s and Spencer~\cite{BES}, with sharp quantitative developments
in \cite{BucicSudakov}. The zero-sum theory was developed by
Caro~\cite{CaroDelta}, Bialostocki and Dierker~\cite{BDMultiple},
and Caro and Provstgaard~\cite{CaroProvstgaard}.

\begin{lemma}\label{thm:balanced-packing}
For every power of two $q\ge2$ and every integer $m\ge2$, put
$g=\gcd(m,q)$. Then
\[
 R(mK_{q,q},\Z_m)=2qm+\frac{m}{g}-1 .
\]
\end{lemma}
\begin{proof}
Put $k=m/g$. The upper bound follows from Caro's transfer inequality
\cite[Theorem~8]{CaroDelta}, which states that
\[
R(tH,\Z_{Mk})\le tR(H,\Z_M)+k-1
\]
for integers $t\ge k\ge2$ with $k\mid t$, and $M\ge1$ with
$M\mid e(H)$. We use the convention $R(H,\Z_1)=|V(H)|$ for the
trivial group, as allowed in that theorem. For $g\ge2$,
Lemma~\ref{thm:balanced-parity} gives $R(K_{q,q},\Z_g)=2q$;
for $g=1$, the same equality follows from the convention.

If $k=1$, then $g=m$. Partition $2qm$ vertices into $m$ sets of
size $2q$ and choose a zero-sum $K_{q,q}$ modulo $m$ on each set.
Their union proves the required upper bound. If $k\ge2$, apply
Caro's inequality with $H=K_{q,q}$, $t=m$ and $M=g$. Since
$k\mid m$ and $g\mid q^2=e(K_{q,q})$, we obtain
\[
R(mK_{q,q},\Z_m)
\le mR(K_{q,q},\Z_g)+\frac{m}{g}-1
=2qm+\frac{m}{g}-1.
\]

For the lower bound, we adapt the degree-count construction from
the proof of~\cite[Theorem~5]{CaroProvstgaard}. If $k=1$, the
vertex-count bound gives $R(mK_{q,q},\Z_m)\ge2qm$.
If $k\ge2$, take $2qm+k-2$ vertices, partition them into sets $A$
and $B$ of sizes $k-1$ and $2qm-1$, and label each edge $e$ by
$|e\cap A|$ modulo $m$. A copy of $mK_{q,q}$ meets $A$ in $h$
vertices for some $1\le h\le k-1$. Every vertex of the packing has
degree $q$, so the packing has total weight $hq$. Since $q$ has additive order
$k$ modulo $m$, this weight is nonzero. Hence
$R(mK_{q,q},\Z_m)\ge2qm+k-1$.
\end{proof}

The same packing argument gives
\[
 R(mK_{n,n},\Z_m)=2nm+\frac{m}{\gcd(m,n)}-1
\]
whenever $n\ge1$ and $\gcd(m,n)$ is a power of two, including $1$.
If $g=\gcd(m,n)\ge2$, Lemma~\ref{thm:rectangular-cuts}
with $q=g$, $k=n/g$ and $s=n$ gives a zero-sum $K_{n,n}$ modulo
$g$ on every $2n$-set. For $g=1$ one can use any balanced cut.
The preceding upper-bound argument and the same degree count then
apply with $n$ in place of $q$.

To compare the packing formula with Conjecture~\ref{conj:defect},
we recall the definition of the defect.

\paragraph{The defect.}
Let $G$ be a graph on $n$ vertices and let $t\ge2$. The number
$\defect_t(G)$ is the least nonnegative integer $d$ such that every
$\Z_{t}$-labeling of $K_{2n}$ has at least one of these properties:
\begin{enumerate}
\item every copy of $G$ has a vertex-disjoint copy of the same weight;
\item some set of at most $n+d$ vertices contains two copies of $G$
with different weights.
\end{enumerate}

When disjoint equal-weight partners are unavailable, the defect
measures how many extra vertices are needed to witness different
weights. Related equal-weight conditions occur in the uniformity
spaces studied by Caro and Yuster~\cite{CaroYusterUniformity}.
Caro~\cite{CaroDelta} proved that one additional vertex always
suffices. We include the short proof for completeness.

\begin{lemma}[Caro~\cite{CaroDelta}]\label{lem:defect-binary}
For every graph $G$ on $n\ge1$ vertices and every $t\ge2$,
$\defect_t(G)\le1$.
\end{lemma}
\begin{proof}
Consider a labeling of $K_{2n}$, and suppose that no set of at most
$n+1$ vertices contains two copies of $G$ with different weights.
All copies on a given $n$-set $S$ then have a common weight $c(S)$.
If $S$ and $S'$ differ by one vertex, their union has $n+1$ vertices,
so $c(S)=c(S')$. Replacing vertices one at a time connects any two
$n$-sets, and therefore every copy of $G$ has the same weight.
The complementary $n$-set supplies a disjoint equal-weight copy
for each copy of $G$. This is the first alternative in
the definition of $\defect_t(G)$, proving the assertion.
\end{proof}

Let $q\ge2$ be a power of two, let $m\ge2$, and put $g=\gcd(m,q)$.
The defect values below follow from the degree criterion of Caro
and Provstgaard~\cite[Theorem~4]{CaroProvstgaard}; we give a direct
verification. If $m\mid q$, then
every $2q$-set has a balanced cut of
weight zero by Lemma~\ref{thm:balanced-parity}. If all cuts on
each $2q$-set have equal weight, every cut therefore has weight
zero and has an equal-weight partner on its complement. This proves
$\defect_m(K_{q,q})=0$. If $m\nmid q$, label the edges incident to
one vertex of $K_{4q}$ by $1$ and all other edges by $0$. A copy of
$K_{q,q}$ containing that vertex has weight $q\ne0$ in $\Z_m$, and
every copy avoiding it has weight zero. Thus all copies on a fixed
$2q$-set have equal weight, while a copy through the marked vertex
has no disjoint equal-weight partner. Lemma~\ref{lem:defect-binary}
now gives $\defect_m(K_{q,q})=1$.

It follows that the original defect formula gives the correct value
if $m\mid q$ or $\gcd(m,q)=1$. In the remaining cases
$1<g<m$, it overestimates the packing threshold by exactly
$m-m/g$.

In particular, for every power of two $q\ge2$, taking $m=2q$ gives
\[
R(2qK_{q,q},\Z_{2q})=4q^2+1,
\qquad \defect_{2q}(K_{q,q})=1.
\]
The conjectured value is $4q^2+2q-1$, so the exact gap is $2q-2>0$.
As explained in Remark~\ref{rem:earlier-work}, the failure of the
conjecture for this family already follows from the earlier results
of Caro and Caro--Provstgaard. For $q=4$, the exact value
$R(8K_{4,4},\Z_8)=65$ sharpens their upper bound $67$; the
conjectured value is $71$.

\subsection{Sufficiently long cycles}
\label{sec:long-cycles}

We prove part~(iii) of Theorem~\ref{thm:main-cuts}. The argument uses
the four-vertex switches and the vertex-label representation from
\cite[Section~4]{ChiHeCycles}. We reserve disjoint switches before
choosing the cycle vertices, and use the Erd\H{o}s--Ginzburg--Ziv theorem
to control both the vertex sum and the total error of the selected paths.
We first record two elementary lemmas.

\begin{lemma}\label{lem:multiple-egz}
Let $A$ be a cyclic group of order $r$, and let $t$ be a positive
multiple of $r$. Every sequence of $t+r-1$ elements of $A$ contains
$t$ terms with sum zero.
\end{lemma}
\begin{proof}
For $r=1$ every sum is zero. For $r\ge2$, repeatedly apply the
Erd\H{o}s--Ginzburg--Ziv theorem \cite{EGZ} to remove $r$ terms of
sum zero. After $j<t/r$ removals, at least
$t+r-1-jr\ge2r-1$ terms remain, so another removal is possible.
The union of the $t/r$ chosen blocks has the required size and sum.
\end{proof}

\begin{lemma}\label{lem:cycle-potentials}
Let $A$ be an abelian group and let $U$ have at least four vertices.
Write $K_U$ for the complete graph on $U$.
Suppose $w:E(K_U)\to A$ satisfies
\[
 w(xz)+w(yt)=w(xy)+w(zt)
\]
for all ordered quadruples of distinct vertices $x,y,z,t$.
Then there are $\alpha_u\in A$ for $u\in U$ and $c\in A$ such that
$w(uv)=\alpha_u+\alpha_v+c$ for every edge.
\end{lemma}
\begin{proof}
Fix $v\in U$. For each $u\ne v$, the difference
$w(ut)-w(vt)$ is independent of $t\notin\{u,v\}$, by the displayed
identity. Denote it by $\alpha_u$ and put $\alpha_v=0$.
For distinct $u,t\ne v$, symmetry gives
\[
 \alpha_u+w(vt)=w(ut)=\alpha_t+w(vu).
\]
Thus $w(vu)-\alpha_u$ has a common value $c$, and the asserted formula
follows. This is the vertex-label argument of
\cite[Lemma~4.2]{ChiHeCycles}, valid in any abelian group.
\end{proof}

\begin{lemma}\label{lem:long-cycles}
Let $m\ge2$, put $D=m/\gcd(m,2)$, and let $\ell$ be a positive
multiple of $m$ with $\ell\ge100m^3$. Then
\[
R(C_\ell,\Z_m)=\ell+D-1.
\]
\end{lemma}
\begin{proof}
\emph{The lower bound.}
We give a degree-count proof of the lower bound in
\cite[Theorem~1.3]{ChiHeCycles}. If $D=1$, the vertex count gives
$R(C_\ell,\Z_m)\ge\ell$. Suppose $D\ge2$. In $K_{\ell+D-2}$ mark
$D-1$ vertices and label each edge by its number of marked endpoints.
Every $\ell$-cycle contains $h$ marked vertices for some
$1\le h\le D-1$, since only $\ell-1$ vertices are unmarked.
Its weight is $2h$, which is nonzero in $\Z_m$ because the element
$2$ has order $D$. Therefore $R(C_\ell,\Z_m)\ge\ell+D-1$.

\smallskip
\noindent\emph{Switches and a quotient.}
Fix an arbitrary labeling $w:E(K_N)\to\Z_m$, where $N=\ell+D-1$.
For distinct vertices $x_1,x_2,x_3,x_4$, define
\[
 \Delta(x_1,x_2,x_3,x_4)
 =w(x_1x_3)+w(x_2x_4)-w(x_1x_2)-w(x_3x_4).
\]
Replacing the path $x_1x_2x_3x_4$ by $x_1x_3x_2x_4$ changes its
weight by $\Delta$ and keeps its endpoints fixed.

Process the nonzero residues $d\in\Z_m$ in any fixed order.
Initially no vertices are reserved. When processing $d$, if the
unreserved vertices contain $3m$ pairwise vertex-disjoint ordered
quadruples with $\Delta=d$, reserve one such family, denoted
$\mathcal F_d$, and record $d$ in a set $B$. Otherwise make no reservation.
Write $P$ for all reserved vertices. The families are mutually
vertex-disjoint, and
\[
 |P|\le 4(3m)(m-1)\le12m^2.
\]
Let $H$ be the subgroup of $\Z_m$ generated by $B$.

Every residue $d\notin H$ was skipped. At the time it was processed,
the unreserved set contained $V(K_N)\setminus P$, so this latter set
also contains no $3m$ pairwise vertex-disjoint quadruples of difference
$d$. Hence, among the vertices outside $P$, a maximal disjoint family of quadruples of difference
$d$ has fewer than $3m$ members. Let $X_d$ be its vertices.
Maximality means that $X_d$ meets every such quadruple.
Put
\[
 X=P\cup\bigcup_{d\notin H}X_d,\qquad U=V(K_N)\setminus X.
\]
Then
\begin{equation}\label{eq:cycle-exceptional}
 |X|\le24m^2.
\end{equation}
Every ordered quadruple in $U$ has difference in $H$.

Let $Q=\Z_m/H$ and write $\overline w$ for the image of $w$ in $Q$.
Since $|U|\ge N-24m^2>4$, Lemma~\ref{lem:cycle-potentials} gives
\[
 \overline w(uv)=\alpha_u+\alpha_v+c\qquad(u,v\in U,\ u\ne v).
\]
For each $x\in X$, pigeonhole the values
$\overline w(xu)-\alpha_u-c$ as $u$ ranges over $U$.
Choose a value $\alpha_x$ attained on a set $U_x\subseteq U$ with
\begin{equation}\label{eq:cycle-pool}
 |U_x|\ge |U|/|Q|\ge |U|/m.
\end{equation}
Thus the same formula holds on $xu$ whenever $u\in U_x$.
The choices for different $x$ are independent.

\smallskip
\noindent\emph{Choosing the vertices and the switches.}

The subgroup $2Q=\{2a:a\in Q\}$ is cyclic of order
\[
 e=\frac{|Q|}{\gcd(|Q|,2)}.
\]
The natural epimorphism $\Z_m\to Q$ maps $2\Z_m$ onto $2Q$,
so $e\mid D$. Also $e\mid\ell$, since $|Q|\mid m$ and $m\mid\ell$. Apply
Lemma~\ref{lem:multiple-egz} to the $N$ elements $2\alpha_v$ of $2Q$.
Since $N\ge\ell+e-1$, there is a set $S$ of $\ell$ vertices with
\begin{equation}\label{eq:cycle-vertex-sum}
 \sum_{v\in S}2\alpha_v=0 \qquad\text{in }Q.
\end{equation}
Exactly $D-1$ vertices were omitted.

Define the residual label on any edge by
\[
 \rho(uv)=\overline w(uv)-\alpha_u-\alpha_v-c\in Q.
\]
It vanishes on edges within $U$ and on $xu$ with $u\in U_x$.
For each reserved quadruple $F=(x_1,x_2,x_3,x_4)$, define its path error
\[
 \varepsilon(F)=\rho(x_1x_2)+\rho(x_2x_3)+\rho(x_3x_4)\in Q.
\]

For every $d\in B$, at least
$3m-(D-1)\ge2m+1$ members of $\mathcal F_d$ have all their vertices
in $S$. This follows because the quadruples in a family are disjoint.
Since $|Q|\mid m$ and $m+|Q|-1\le2m-1$,
Lemma~\ref{lem:multiple-egz} selects exactly $m$ such members whose
errors sum to zero. Denote the selected family by
$\mathcal E_d$. Thus
\begin{equation}\label{eq:cycle-errors}
 \sum_{F\in\mathcal E_d}\varepsilon(F)=0\qquad(d\in B).
\end{equation}
The $m$ selected paths for each $d$ will provide independent switches.

\smallskip
\noindent\emph{Assembling the cycle.}

For each selected path $x_1x_2x_3x_4$, choose unused vertices
$u\in U_{x_1}\cap S$ and $v\in U_{x_4}\cap S$ and form the longer
path $ux_1x_2x_3x_4v$. For every vertex $x\in X\cap S$ outside the
selected quadruples, choose two unused vertices $u,v\in U_x\cap S$
and form the path $uxv$. Require all the chosen vertices of $U$
to be distinct.

These choices are possible greedily. They need at most $2|X|$ vertices
of $U$. By \eqref{eq:cycle-exceptional}, \eqref{eq:cycle-pool}, and $N\ge100m^3$,
each pool initially satisfies
\[
 |U_x\cap S|
 \ge\frac{N-24m^2}{m}-(D-1)
 \ge100m^2-25m+1
 >48m^2+2
 \ge2|X|+2.
\]

All resulting paths are disjoint and have endpoints in $U$.
Treat every remaining vertex of $S\cap U$ as a one-vertex path.
There are at least $\ell-3|X|\ge100m^3-72m^2>3$ such vertices.
Join all the paths in any cyclic order to obtain a simple cycle $C$
on $S$. Every joining edge has both endpoints in $U$.

All its residual edge labels vanish except those inside selected
quadruples. Their total is zero by \eqref{eq:cycle-errors}. Consequently,
\[
 \overline w(C)
 =\sum_{v\in S}2\alpha_v+\ell c+
   \sum_{d\in B}\sum_{F\in\mathcal E_d}\varepsilon(F)
 =0.
\]
Here \eqref{eq:cycle-vertex-sum} gives the first zero, and $m\mid\ell$ gives
$\ell c=0$ in $Q$. Therefore $w(C)\in H$.

\smallskip
\noindent\emph{Adjusting the weight.}

Since $B$ generates $H$, there are coefficients
$a_d\in\{0,1,\ldots,m-1\}$ such that
\[
 -w(C)=\sum_{d\in B}a_d d\qquad\text{in }\Z_m.
\]
For each $d\in B$, perform the interior swap on exactly $a_d$ of
the $m$ paths selected from $\mathcal F_d$. Each swap preserves its
endpoints, and the paths are disjoint, so the operations leave a
single cycle on the same $\ell$ vertices. Its new weight is zero.
This proves the required upper bound.
\end{proof}

\Needspace{7\baselineskip}
\begin{proof}[Proof of Theorem~\ref{thm:main-cuts}]
Part~(i), including the threshold $q-2$, follows from
Lemma~\ref{thm:rectangular-cuts} with $\Gamma=\Z_q$.
Part~(ii) is Lemma~\ref{thm:balanced-packing}, and
part~(iii) is Lemma~\ref{lem:long-cycles}.
\end{proof}

\Needspace{13\baselineskip}
\section*{Declaration of generative AI and AI-assisted technologies}
The authors used a large language model during the preparation of
this manuscript for assistance with mathematical exploration,
checking arguments, organization and language editing. The authors
are responsible for the mathematical statements, proofs, references
and final text.

\end{document}